\documentclass[draft,reqno,11pt,oneside]{article} 

\usepackage{amssymb}
\usepackage{english}
\usepackage{amsthm}
\usepackage{amssymb}
\usepackage{longtable} 
\usepackage{amsmath}
\usepackage{ifthen} 
\usepackage{upref}
\usepackage{bbold} 
\usepackage{latexsym} 
\usepackage[all]{xy}
\usepackage{fancyhdr} 
\usepackage{makeidx} 
\newfont{\chen}{cmex10 at 9pt} \newfont{\chee}{cmex10 at 10pt} 
\newfont{\cheu}{cmex10 at 11pt} 
\makeindex 
\begin{document} 
\newtheorem{lemma}{Lemma}[section]
\newtheorem*{eldemo}{Sketch of proof}
\newtheorem{satz}[lemma]{Theorem}
\newtheorem{prop}[lemma]{Proposition}
\newtheorem{defi}[lemma]{Definition}
\newtheorem{bei}[lemma]{Example} 
\newtheorem{verm}[lemma]{Conjecture} 
\newtheorem{kor}[lemma]{Corollary} 
\renewcommand{\proofname}{Proof} 
\newtheorem{bez}[lemma]{Notation} 
\newtheorem{bem}[lemma]{Remark}
\newtheorem{fall}[lemma]{Case}
\newtheorem*{example}{Example} 
\newtheorem{defsatz}[lemma]{Definition and Theorem} 
%
\newtheorem{lemmaap}{Lemma}[section]
\newtheorem{satzap}[lemmaap]{Theorem}
\newtheorem{propap}[lemmaap]{Proposition}
\newtheorem{defiap}[lemmaap]{Definition}
\newtheorem{beiap}[lemmaap]{Example}  
\newtheorem{korap}[lemmaap]{Corollary}  
\newtheorem{bezap}[lemmaap]{Notation} 
\newtheorem{bemap}[lemmaap]{Remark}
\def \Q{{\mathbb{Q}}}
\def \Ar{{\mathcal{A}}}
\def \G{{\mathcal{G}}} 
\def \pu{{\mbox{.}}} 
\def \ps{{\Phi^s}} 
\def \LUC{{\mathrm{LUC}}} 
\def \RUC{{\mathrm{RUC}(\G)}}
\def \RSp{{\widetilde{R}}^{(p)}} 
\def \RSq{{\widetilde{R}}^{(q)}} 
\def \wwp{{\mathfrak{P}}} 
\def \phii{{\overline{\Phi}}} 
\def \Mun{{\mathcal{M}}}
\def \Ce{{\mathrm{C}}}
\def \Cede{{\mathbb{C}}} 
\def \lpp{{L_p}} 
\def \Te{{\mathcal{T}}}
\def \sis{{\mathcal{S}}} 
\def \th{{\widetilde{\Theta}}} 
\def \mc{{\mathcal{MC}}} 
\def \ml{{\mathcal{MLUC}}} 
\def \sus{{\subseteq}} 
\def \mr{{\mathcal{MRUC}}} 
\def \N{{\mathbb{N}}}
\def \al{{\alpha}} 
\def \a{{\alpha}} 
\def \suppp{{\mathrm{supp}}} 
\def \cb{{\mathrm{cb}}} 
\def \A{{\mathcal{A}}}
\def \llu{{\mathrm{LUC}(\G)}} 
\def \B{{\mathcal{B}}}
\def \kapp{{\mathfrak{k}}} 
\def \ru{{{\RUC}^*}} 
\def \L{{\mathfrak{L}}}
\def \EL{{\mathcal{L}}}
\def \gip{{\Gamma_{l, p}}} 
\def \UC{{{\mathrm{UC}}(\G)}} 
\def \gup{{{\widetilde{\Gamma_{l, p}}}}} 
\def \gurp{{\Gamma_{r, p}}} 
\def \gurrp{{{\widetilde{\Gamma_{r, p}}}}} 
\def \M{{\rm{M}}}
\def \en{{\mathcal{N}}}
\def \gr{{\widetilde{\Gamma_r}}} 
\def \gre{{\Gamma_r}} 
\def \em{{\mathcal{M}}} 
\def \emli{{\mathcal{ML}_\infty}} 
\def \ccc{{\mathcal{C}}}
\def \Ha{{\mathfrak{H}}}
\def \l1{{L_1(\G)}} 
\def \li{{L_{\infty}(\G)}} 
\def \H{{\mathcal{H}}}
\def \Qu{{\mathfrak{Q}}}
\def \cbgl{~{\stackrel{\mathrm{cb}}{=}}~} 
\def \FS{{\mathcal{FS}}} 
\def \gammas{{\widetilde{\gamma}}} 
\def \C{{\mathcal{C}}}
\def \mr{{\mathcal{MRUC}}}
\def \Be{{\mathfrak{B}}}
\def \bar{{~|~}} 
\def \Ka{{\mathcal{K}}}
\def \rr{{\mathcal{R}}} 
\def \ceb{{\C\B(\B(L_2(\G)))}} 
\def \id{{{\mathrm{id}}}} 
\def \cep{{\C\B(\B(L_p(\G)))}} 
\def \ep{{\B(\B(L_p(\G)))}}
\def \Is{{\widetilde{I}}} 
\def \conv{{\mathrm{conv}}} 
\def \convs{{\widetilde{\mathrm{conv}}}} 
\def \convsp{{\widetilde{\mathrm{conv_p}}}} 
\def \es{{\B(\B(L_2(\G)))}}
\def \go{{\gip}} 
\def \tt{{\overline{\otimes}}} 
\def \ti{{\stackrel{\vee}{\otimes}}} 
\def \oo{{{\overline{\otimes}}}} 
\def \tp{{\widehat{\otimes}}} 
\def \otp{{\otimes}} 
\def \zet{{Z_t(\lu)}} 
\makeatletter 
\newcommand{\essu}
{\mathop{\operator@font ess\mbox{-}sup}} 
\newcommand{\essi}
{\mathop{\operator@font ess\mbox{-}inf}} 
\newcommand{\lisu}
{\mathop{\operator@font lim\,ess\mbox{-}sup}} 
\newcommand{\liin}
{\mathop{\operator@font lim\,ess\mbox{-}inf}} 
\makeatother 
\def \ga{{\Gamma_p}}
\def \gam{{\gamma_p}}
\def \gaw{{\widetilde{\ga}}}
\def \lu{{\LUC(\G)^*}}
\def \sa{{\overline{\Gamma_p(\M(\G))}^{w*}}}
\def \Ball{{\mathrm{Ball}}} 
\def \gurr{{\widetilde{\Gamma_r}}} 
\def \gur{{\Gamma_r}} 
\def \Ker{{\it{KERN}}}
\def \F{{\mathfrak{F}}}
\def \ad{{\mathrm{ad}}} 
\def \Ceee{{\mathbb{C}}}
\def \Cee{{\mathrm{C}}}
\def \imink{{\iota_{\mathrm{minK}}}} 
\def \imin{{\iota_{\mathrm{min}}}} 
\def \Bild{{\it{BILD}}}
\def \Sz{{\mathcal{S}}}
\def \The{{\mathcal{T}}}
\def \cebe{{\C\B(\B(\Sz_2))}}
\def \cebes{{\C\B(\B(\Sz_2(L_2(G))))}} 
\def \pr{{\Sz_2 \otimes_h \B(\Sz_2) \otimes_h \Sz_2}} 
\def \su{{\rm{sup}}_{t \in \G}}
\def \mi{{\Gamma_2(\M(\G))}}
\def \Cb{{\mathcal{C}}}
\newcommand{\ens}{{\mathcal{N}}(L_p(\G))}
\newcommand{\Tee}{{\mathcal{T}}(\mathcal{H})}
\newcommand{\schutz}[1]{#1} 
\def \ma{{\Gamma_2(\lu)}}
\def \Hi{{\mathcal{H}}}
\def \bs{{\B^{\sigma}(\B(\Ha))}}
\def \bsa{{\C\B(\B(\H))}}
\def \bsi{{\B^s(\B(\Ha))}}
\def \gu{{\widetilde{\Gamma_l}}}
\def \EF{{\widetilde{F}}} 
\def \gi{{\Gamma_l}} 
\def \als{{\widetilde{\alpha}}} 
\def \betas{{\widetilde{\beta}}} 
\def \gir{{\Gamma_r}} 
\def \muss{{\widetilde{\mu}}} 
\def \ge{{\Gamma}}
\def \gem{{\Gamma(\M(\G))}}
\def \MM{{\mathbf{M}}} 
\def \un{{\ell_{\infty}^*(\G)}}
\def \R{{\mathcal{R}}}
\def \cebr{{\C\B_{\R(\G)}(\B(L_2(\G)))}}
\def \Chi{{\chi}} 
\def \cebre{{\C\B_{\R(\G)}(\B(\ell_2(\G)))}}
\def \cebrs{{\C\B_{\R(\G)^{'}}^{\sigma}(\B(L_2(\G)))}}
\def \cebra{{\C\B_{\R(\G)}^{\sigma}(\B(\ell_2(\G)))}}
\def \U{{\mathfrak{U}}}
\def \K{\Ka} 
\def \cebru{{\C\B_R(\B(L_2(\G)))}}
\def \lin{{\mathrm{lin}}} 
\def \cebri{{\C\B_{R}^{\sigma}(\B(\H))}}
\def \Tee{{\The(\H)}}
\def \gurrn{{\gurr^0}} 
\def \gin{{\gi^0}} 
\def \gun{{\gu^0}} 
\newcommand{\dogl}{
\setlength{\unitlength}{0.7ex}
\linethickness {0.1ex}
~\begin{picture} (2.90 , 2)
\put (0 , 0.5) {\circle* {0.3}}
\put (0 , 1.2) {\circle* {0.3}}
\put (0.7 , 0.5) {\line (1,0) {2.1}}
\put (0.7 , 1.2) {\line (1,0) {2.1}}
\end{picture}~}
%
%
\newlength{\checklength}
\newcommand{\fns}{\footnotesize}
\newcommand{\widecheck}[1]{
\settowidth{\checklength}{$#1$}
\linethickness {0.1ex}
\setlength{\unitlength}{\checklength}
\stackrel{
\ifthenelse{ \lengthtest{ \checklength > 6em } }
{
\begin{picture}(1,0.04)
\qbezier(0,0.04)(0.3,0.03)(0.5,0.003)
\qbezier(0.5,0.003)(0.7,0.03)(1,0.04)
\qbezier(0,0.04)(0.3,0.03)(0.5,-0)
\qbezier(0.5,-0)(0.7,0.03)(1,0.04)
\qbezier(0,0.04)(0.3,0.03)(0.5,-0.005)
\qbezier(0.5,-0.005)(0.7,0.03)(1,0.04)
\end{picture}
}
{
\ifthenelse{ \lengthtest{ \checklength > 4em } }
{
\begin{picture}(1,0.05)
\qbezier(0,0.05)(0.3,0.03)(0.5,-0.02)
\qbezier(0.5,-0.02)(0.7,0.03)(1,0.05)
\qbezier(0,0.05)(0.3,0.03)(0.5,-0.015)
\qbezier(0.5,-0.015)(0.7,0.03)(1,0.05)
\qbezier(0,0.05)(0.3,0.03)(0.5,-0.01)
\qbezier(0.5,-0.01)(0.7,0.03)(1,0.05)
\end{picture}
}
{
\ifthenelse{\lengthtest{\checklength > 3em}}
{
\begin{picture}(1,0.075)
\qbezier(0,0.075)(0.3,0.05)(0.5,-0.03)
\qbezier(0.5,-0.03)(0.7,0.05)(1,0.075)
\qbezier(0,0.075)(0.3,0.05)(0.5,-0.02)
\qbezier(0.5,-0.02)(0.7,0.05)(1,0.075)
\end{picture}
}
{
\ifthenelse{\lengthtest{\checklength > 2em}}
{
\begin{picture}(1,0.1)
\qbezier(0,0.1)(0.3,0.07)(0.5,-0.04)
\qbezier(0.5,-0.04)(0.7,0.07)(1,0.1)
\qbezier(0,0.1)(0.3,0.07)(0.5,-0.02)
\qbezier(0.5,-0.02)(0.7,0.07)(1,0.1)
\end{picture}
}
{
\ifthenelse{\lengthtest{\checklength > 1em}}
{
\begin{picture}(1,0.12)
\qbezier(0,0.12)(0.3,0.1)(0.5,-0.06)
\qbezier(0.5,-0.06)(0.7,0.1)(1,0.12)
\qbezier(0,0.12)(0.3,0.1)(0.5,-0.025)
\qbezier(0.5,-0.025)(0.7,0.1)(1,0.12)
\end{picture}
}
{
\begin{picture}(1,0.15)
\qbezier(0,0.15)(0.3,0.12)(0.5,-0.1)
\qbezier(0.5,-0.1)(0.7,0.12)(1,0.15)
\qbezier(0,0.15)(0.3,0.12)(0.5,-0.04)
\qbezier(0.5,-0.04)(0.7,0.12)(1,0.15)
\end{picture}
}
}
}
}
}
}
{#1}
}
\newcommand{\swidecheck}[1]{
\settowidth{\checklength}{$_{#1}$}
\linethickness {0.1ex}
\setlength{\unitlength}{0.95\checklength}
\stackrel{
\ifthenelse{ \lengthtest{ \checklength > 6em } }
{
\begin{picture}(1,0.04)
\qbezier(0,0.04)(0.3,0.03)(0.5,0.003)
\qbezier(0.5,0.003)(0.7,0.03)(1,0.04)
\qbezier(0,0.04)(0.3,0.03)(0.5,-0)
\qbezier(0.5,-0)(0.7,0.03)(1,0.04)
\qbezier(0,0.04)(0.3,0.03)(0.5,-0.005)
\qbezier(0.5,-0.005)(0.7,0.03)(1,0.04)
\end{picture}
}
{
\ifthenelse{ \lengthtest{ \checklength > 4em } }
{
\begin{picture}(1,0.05)
\qbezier(0,0.05)(0.3,0.03)(0.5,-0.02)
\qbezier(0.5,-0.02)(0.7,0.03)(1,0.05)
\qbezier(0,0.05)(0.3,0.03)(0.5,-0.015)
\qbezier(0.5,-0.015)(0.7,0.03)(1,0.05)
\qbezier(0,0.05)(0.3,0.03)(0.5,-0.01)
\qbezier(0.5,-0.01)(0.7,0.03)(1,0.05)
\end{picture}
}
{
\ifthenelse{\lengthtest{\checklength > 3em}}
{
\begin{picture}(1,0.075)
\qbezier(0,0.075)(0.3,0.05)(0.5,-0.03)
\qbezier(0.5,-0.03)(0.7,0.05)(1,0.075)
\qbezier(0,0.075)(0.3,0.05)(0.5,-0.02)
\qbezier(0.5,-0.02)(0.7,0.05)(1,0.075)
\end{picture}
}
{
\ifthenelse{\lengthtest{\checklength > 2em}}
{
\begin{picture}(1,0.1)
\qbezier(0,0.1)(0.3,0.07)(0.5,-0.04)
\qbezier(0.5,-0.04)(0.7,0.07)(1,0.1)
\qbezier(0,0.1)(0.3,0.07)(0.5,-0.02)
\qbezier(0.5,-0.02)(0.7,0.07)(1,0.1)
\end{picture}
}
{
\ifthenelse{\lengthtest{\checklength > 1em}}
{
\begin{picture}(1,0.12)
\qbezier(0,0.12)(0.3,0.12)(0.5,-0.06)
\qbezier(0.5,-0.06)(0.7,0.12)(1,0.12)
\qbezier(0,0.12)(0.3,0.12)(0.5,-0.025)
\qbezier(0.5,-0.025)(0.7,0.12)(1,0.12)
\end{picture}
}
{
\begin{picture}(1,0.15)
\qbezier(0,0.15)(0.3,0.12)(0.5,-0.1)
\qbezier(0.5,-0.1)(0.7,0.12)(1,0.15)
\qbezier(0,0.15)(0.3,0.12)(0.5,-0.04)
\qbezier(0.5,-0.04)(0.7,0.12)(1,0.15)
\end{picture}
}
}
}
}
}
}
{#1}
} 
\def \des{{\sqrt{x_j} \otimes M_{T_{\widecheck{b}}} \otimes \sqrt{x_j}}} 
\def \res{{(e_j)^{\frac{1}{p}} \otimes T_{\widecheck{b}} \otimes (e_j)^{\frac{1}{q}}}} 
\def \rem{{(e_j)^{\frac{1}{2}} \otimes M_{\widecheck{b}} \otimes (e_j)^{\frac{1}{2}}}} 
%
%
\title{A Unified Approach to the Topological Centre 
Problem \\ 
for Certain Banach Algebras Arising in Abstract Harmonic Analysis} 
\author{Matthias Neufang}
\date{} 
\maketitle 
\begin{abstract} 
Let $\G$ be a locally compact group. Consider the Banach algebra $\l1^{**}$, equipped with the 
first Arens multiplication, as well as the algebra $\lu$, the dual of the space of 
bounded left uniformly continuous functions on $\G$, whose product extends the convolution in the measure algebra $\M(\G)$. 
We\footnote{ 
\par 
The author is currently a PIMS Postdoctoral Fellow at the 
University of Alberta, Edmonton, where this work was accomplished. The support of PIMS is gratefully acknowledged.} 
present (for the most interesting case of a non-compact group) 
completely different -- in particular, \textit{direct} -- proofs 
and even obtain \textit{sharpened} versions of the results, first proved by Lau--Losert in 
\cite{lalo} and Lau in \cite{lau}, that the topological centres of the latter algebras precisely are 
$\l1$ and $\M(\G)$, respectively. The 
special interest of our new approach lies in the fact that it shows a 
fairly general pattern of solving the topological centre problem for various kinds of Banach algebras; in particular, 
it 
avoids the use of any measure theoretical techniques. At the same time, 
deriving both results in perfect parallelity, 
our method 
reveals the nature of 
their 
close relation. 
\end{abstract}

\section{Introduction} 
In this note, we wish to present a 
new and, for the first time, \textit{unified} approach to two theorems which 
may be considered as the fundamental results, known by now, 
concerning the topological centre problem 
for concrete Banach algebras studied in Abstract Harmonic Analysis. 
Namely, for a period of nearly 15 years (beginning in the seventies with \cite{zap}), 
research in the topological centre question 
was centered around the Banach algebras 
$L_1(\G)^{**}$, endowed with the first Arens product, and its quotient algebra $\lu$, where $\G$ denotes a 
locally compact group. -- Here, we write $\LUC (\G)$ for the space of complex-valued bounded left uniformly continuous functions 
on $\G$; the definition of the product on the dual $\lu$ will be briefly recalled below. 
The questions were eventually answered in full generality by the decisive work of Lau \cite{lau} and Lau--Losert \cite{lalo}. 
In the present note, we will derive both these results (see Theorems \ref{hau} and \ref{hau10} below) by organizing the 
arguments in a parallel fashion, 
which at the same time 
not only will yield direct proofs, but even sharpenings of the statements. 
\begin{satz} \label{hau}
The topological centre $Z_t(L_1(\G)^{**})$ of the Banach algebra $L_1(\G)^{**}$ is precisely $L_1(\G)$. 
\end{satz}
\begin{proof}
This 
is the main result, Thm.\ 1, in \cite{lalo}. 
\end{proof} 
We shall restrict ourselves to the most interesting situation where $\G$ is non-compact. In case $\G$ is compact, 
the above assertion may be attained by a very short argument (\cite{lalo}, Thm.\ 1; also cf.\ \cite{ipu}, Thm.\ 3.3 (vi)). 
\par 
The second theorem even is only of interest in the non-compact case -- for in case $\G$ is compact, 
$\lu$ equals the measure algebra $\M(\G)$, and the assertion is immediately veryfied. 
\begin{satz} \label{hau10}
The topological centre $Z_t(\lu)$ of the Banach algebra $\lu$ is $\M (\G)$. 
\end{satz}
\begin{proof}
This is the main result, Thm.\ 1, in \cite{lau}. 
\end{proof} 
In the sequel, we will have to distinguish between the two Arens products on the bidual $\l1^{**}$. We 
denote by $\odot$ the first ($=$ left) and by $.$ the second ($=$ right) Arens product, 
and use these symbols as well for the various module operations linking $\l1$, its dual and bidual, as follows. 
Let $m , n \in \l1^{**}$, $h \in \l1^*$, $f , g \in \l1$. Denoting by $*$ the convolution product of functions 
(whenever it is defined), we write: 
\begin{center} 
\begin{tabular}{ccc} 
$\langle h \odot f, g \rangle$ & $:=$ & $\langle h, f * g \rangle$ \\ 
$\langle n \odot h, f \rangle$ & $:=$ & $\langle n, h \odot f \rangle$ \\ 
$\langle m \odot n, h \rangle$ & $:=$ & $\langle m, n \odot h \rangle$ 
\end{tabular} 
\end{center} 
and, following a completely symmetric pattern: 
\begin{center} 
\begin{tabular}{ccc}
$\langle f . h, g \rangle$ & $:=$ & $\langle h, g * f \rangle$ \\ 
$\langle h . m, f \rangle$ & $:=$ & $\langle m, f . h \rangle$ \\ 
$\langle m . n, h \rangle$ & $:=$ & $\langle n, h . m \rangle$. 
\end{tabular} 
\end{center} 
We note that, in particular, 
$h \odot f = \widetilde{f} * h \in \LUC(\G)$; here, 
as usual, we write $\widetilde{f} := \frac{1}{\Delta} \check{f}$, where $\Delta$ denotes the modular function of $\G$ and 
$\check{f} (x) := f(x^{-1})$ for all $x \in \G$. 
We refer to \cite{pal}, 
\S 1.4, for a discussion of basic properties of Arens multiplication in the framework of general Banach algebras. 
\par 
Whenever 
we 
consider $\l1^{**}$ as a Banach algebra, we regard it as 
equipped with the first Arens product. 
We recall that the 
topological centre of $\l1^{**}$ 
is defined to be the set of functionals $m \in \l1^{**}$ which satisfy $m \odot n = m . n$ for all 
$n \in \l1^{**}$. Equivalently, the topological centre consists of all the functionals $m \in \l1^{**}$ such 
that left multiplication by $m$ is $w^*$-$w^*$-continuous on $\l1^{**}$. 
A detailed analysis of topological centres in the general context of biduals of Banach algebras 
can be found, e.g., in \cite{laul}. 
\par 
There is an analogous notion of topological centre for the Banach algebra $\lu$. First, let us recall the natural 
construction of the product 
in the latter space. -- If $n \in \lu$ and $f \in \LUC(\G)$, then it is classical that 
the function $n \cdot f$, defined through $$(n \cdot f)(x) := \langle n , l_x f \rangle \quad (x \in \G)$$
still belongs to $\LUC(\G)$; i.e., $\LUC(\G)$ is left introverted. 
This operation 
gives rise to the 
product on the space $\lu$: 
$$\langle m \cdot n , f \rangle := \langle m , n \cdot f \rangle \quad (m , n \in \lu , f \in \LUC(G))$$ 
under which 
the latter indeed becomes a Banach algebra, and $\LUC(\G)$ becomes a left $\lu$-module with the action introduced above. 
In analogy to the case of $\l1^{**}$, one defines the topological centre of $\lu$ to be the set of 
elements $m \in \lu$ such that left multiplication by $m$ is $w^*$-$w^*$-continuous on $\lu$. 
\par 
We shall also consider the natural (left) module operation 
of $\lu$ on $\li$ given by 
$\langle m \Diamond h, g \rangle := \langle m, h \odot g \rangle$, 
where $m \in \lu$, $h \in \li$, $g \in \l1$. 
At this point, we note that, as is easily seen, 
one has $m \Diamond h = \widetilde{m} \odot h$, 
where $\widetilde{m}$ is an arbitrary Hahn-Banach extension of the functional $m$ to $\li^*$. 
\par 
The main interest of the approach presented here 
consists in the following: 
\begin{itemize}
\item 
We obtain a 
sharpening of the non-trivial inclusion in the 
statements of both Theorem \ref{hau} and \ref{hau10}. Namely, as for Theorem \ref{hau}, we shall see that for an element 
$m \in \li^*$ in order to belong to $\l1$, it is sufficient 
to have $m \odot n = m . n$ only for all $n \in \li^*$ that are Hahn-Banach extensions of functionals in 
${\overline{\delta_\G}}^{w^*} \subseteq \Ball(\lu)$ -- instead of requiring the latter equality for all $n \in \li^*$, 
as does the definition of the topological centre. 
Analogously, in the situation of Theorem \ref{hau10}, we will prove that an element $m \in \lu$ already 
belongs to $\M(\G)$ if left multiplication by $m$ is only 
$w^*$-$w^*$-continuous on ${\overline{\delta_\G}}^{w^*} \subseteq \Ball(\lu)$, instead of demanding this 
continuity on the whole unit ball of $\lu$. 
\item 
The only proof known so far of the topological centre theorem for $\lu$ is indirect (see Thm.\ 1 in \cite{lau}), and all 
proofs given for the corresponding theorem for $\l1^{**}$ either heavily rely on the latter (\cite{laul}, \cite{ghmc}) or are also 
indirect (see \cite{lalo}, Thm.\ 1). Our 
proofs of the two results 
are independent, and in both cases \textit{direct}. -- We remark \textit{en passant} that the proofs of Thm.\ 5.4 and Cor.\ 5.5 
in \cite{laul} are correct only under the additional set-theoretic assumption that the compact covering number of the 
group $\G$ is a non-measurable cardinal, since the argument is based on \textit{ibid.}, Lemma 5.3, which in turn has to be read 
with a similar set-theoretic assumption; we refer to \cite{ich1} and \cite{hune} for a detailed discussion of these and related problems. 
The tool which enables us to overcome precisely those set-theoretic difficulties is provided by Lemma \ref{kmazur} below. 
\item 
The procedure we present shows a perfect analogy 
between the two topological centre theorems. Besides one additional structural result for $\lu$, the prerequisites are 
the same, and the proofs themselves follow completely parallel lines. 
\item 
The method of proof follows a purely Banach algebraic procedure and does not, in particular, rely on any 
measure theoretic argument, 
so that it might be applied 
equally well in other situations. 
At this point, we only mention that in \cite{hune}, we obtained an analogue 
of Lemma \ref{kmazur} below in the ``dual'' setting, i.e., with $\li$ replaced by the group von Neumann algebra 
$\mathit{VN} (\G)$, and with $\kapp (\G)$ replaced by $b (\G)$, the smallest cardinality of an open basis at the neutral element of $\G$. 
For Proposition \ref{ftcsatz} as well, 
a dual version is known at least in the case of amenable groups $\G$ -- here, 
the algebras $\mathrm{UCB}(\widehat{\G})^*$ and $B_\rho (\G)$ take over the place of $\lu$ and $\M(\G)$, respectively. Thus, 
if 
we are able to prove a 
factorization result, for bounded families in $\mathit{VN} (\G)$ of cardinality at most 
$b (\G) \cdot \aleph_0$, corresponding to Lemma \ref{fakt}, then our method of proof 
would immediately yield an affirmative answer 
(for the most interesting case of non-discrete groups $\G$) 
of two longstanding 
conjectures at once -- namely, the topological centre of the bidual of the Fourier algebra $A(\G)$ being just $A(\G)$, i.e., 
$Z_t (A(\G)^{**}) = A(\G)$, and, for amenable groups $\G$, $Z_t (\mathrm{UCB}(\widehat{\G})^*) = B_\rho (\G)$. 
\end{itemize}

\section{Preliminaries} 
For both proofs 
we will use the following two lemmata, which are of interest in their own right. 
\begin{lemma} \label{kmazur}
For an arbitrary locally compact group $\G$, the space $\l1$ enjoys Mazur's property of 
level $\kapp(\G) \cdot \aleph_0$, where $\kapp(\G)$ denotes the 
compact covering number of $\G$ 
(i.e., the least cardinality of 
a compact covering of $\G$). 
-- This means that 
a functional $m \in \l1^{**}$ actually belongs to $\l1$ if it carries 
bounded $w^*$-converging nets of cardinality at most $\kapp(\G) \cdot \aleph_0$ into converging nets. 
\end{lemma}
\begin{proof}
This is Thm.\ 4.4 in \cite{ich1}. 
\end{proof} 
\noindent 
Next we present 
our crucial tool, which is a general factorization theorem for 
bounded families in $\li$. It has already been used (cf.\ \cite{ich2}) 
\begin{itemize} 
\item 
to answer (in the affirmative) a question raised by Hofmeier--Wittstock in \cite{wiho} concerning the automatic 
boundedness of left $\li^*$-module homomorphisms on $\li$ 
\item 
to give an alternative approach to the result on 
automatic $w^*$-$w^*$-continuity of the latter mappings, as first shown by 
Ghahramani--McClure in \cite{ghmc}. 
\end{itemize} 
\begin{lemma} \label{fakt}
Let $\G$ be a locally compact non-compact group with compact covering number $\kapp(\G)$. Let further 
$(h_\al)_{\al \in I} \subseteq \li$ be a bounded family of functions where $|I| \leq \kapp(\G)$. Then there exist 
a family $(\psi_\al)_{\al \in I}$ of functionals in ${\overline{\delta_\G}}^{w^*} \subseteq \Ball(\lu)$ and a 
function $h \in \li$ such that the factorization formula 
$$h_\al = \psi_\al \Diamond h$$
holds for all $\al \in I$. (Moreover, the functionals $\psi_\al$, $\al \in I$, do not depend, except for the 
index set, on the given family $(h_\al)_{\al \in I}$; they are universal in the sense that they 
are obtained intrinsically from the 
group $\G$.) 
\end{lemma}
\begin{proof}
See 
\cite{ich}, Satz 3.6.2 (or \cite{ich2}). 
\end{proof} 
In order to arrange the proof of Thm.\ \ref{hau10} in a completely 
parallel manner, 
we only require 
the following proposition, whose first part is a classical 
structural result about the algebra $\lu$ going back to the 
pioneering work of Curtis--Fig\`{a}-Talamanca (\cite{cufi}, Thm.\ 3.3). We write $\B^\sigma (L_\infty (\G))$ 
for the space of normal (i.e., $w^*$-$w^*$-continuous) operators on $L_\infty (\G)$. 
%
\begin{prop} \label{ftcsatz} 
The mapping 
$$\phi : \lu \longrightarrow \B (L_{\infty}(\G))$$ 
defined through 
$$\phi (m) (h) := m \Diamond h \quad ( m \in \lu, ~h \in \li )$$ 
is an isometric representation of $\lu$ in $\B (\li)$ such that 
\begin{eqnarray} 
\phi (\M(\G)) = \phi (\lu) \cap \B^{\sigma}(\li) . \label{rel} 
\end{eqnarray} 
\end{prop}
\begin{proof}
The first assertions are well-known (see, e.g., \cite{tlau}, Thm.\ 1, together with Lemma 1 and Remark 3; or 
\cite{cufi}, Thm.\ 3.3, for a proof in the unimodular case). 
The relation (\ref{rel}) is seen as follows (cf.\ \cite{ich}, Prop.\ 3.1.1). -- In order to prove the inclusion ``$\subseteq$'', 
let $\mu \in \M(\G)$. We have to show that $\phi (\mu) \in \B(\li)$ is normal. 
Consider a net $(h_\al)_\al \subseteq \Ball(\li)$ such that $h_\al \stackrel{w^*}{\longrightarrow} 0$. Fix $g \in \l1$. 
We claim that 
$$\langle \phi (\mu)(h_\al), g \rangle \longrightarrow 0.$$ 
But we have: 
$$\langle \phi (\mu)(h_\al), g \rangle = \langle \mu \Diamond h_\al , g \rangle = \langle \mu, {\widetilde{g}} * h_\al \rangle 
= \int_\G ({\widetilde{g}} * h_\al)(t) ~d\mu (t) ,$$ 
where (as is readily checked) 
the net 
$({\widetilde{g}} * h_\al)_\al \subseteq \LUC(\G)$ is equicontinuous, bounded 
and converges pointwise to $0$. It thus converges uniformly on compact subsets of $\G$, 
whence we conclude that the above integrals converge to $0$, as desired. 
\par 
Turning to the inclusion ``$\supseteq$'', let us consider 
an arbitrary element $\phi (m)$ of the set on the right side, where $m \in \lu$. 
We denote by $\Cee_0(\G)$ the space of all complex-valued continuous functions on $\G$ 
vanishing at infinity, and by $\Cee_0(\G)^{\perp}$ its annihilator in $\lu$. Then 
we obviously have $\lu = \M(\G) \oplus \Cee_0(\G)^{\perp}$ (the latter is actually 
an $\ell_1$-direct sum -- cf.\ Lemma 1.1 in \cite{gll} --, but we will not need this fact). 
Write $m=\mu + n$ with $\mu \in \M(\G)$ and $n \in \Cee_0(\G)^{\perp}$ according to this decomposition. 
It suffices to show that $\phi (n)=0$. 
\par 
First, using the inclusion proved above, we see that $\phi (n) = \phi (m) - \phi (\mu)$ is normal. Fix $h \in \li$. 
We have $h = \sigma(\li, \l1)-\lim_\al h_\al$ for 
an appropriate net $(h_\al)_\al \subseteq \Cee_0(\G)$. 
Hence we obtain: 
$$\phi (n)(h) = \sigma(\li, \l1)-\lim_\al \phi (n)(h_\al) = 0 ,$$ 
which finishes the proof. -- For the last equality, note that 
for all $\al$ and arbitrary $g \in \l1$: 
$$\langle \phi (n)(h_\al), g \rangle = \langle n, {\widetilde{g}} * h_\al \rangle = 0,$$ 
since ${\widetilde{g}} * h_\al \in L_1(\G) * \Cee_0(\G) = \Cee_0(\G)$. 
\end{proof} 
\begin{bem} 
We only stated the specific properties of $\phi$ for the sake of completeness. In the sequel, we shall just use the fact 
that $\phi$ is bounded and injective, as well as the relation $(\ref{rel})$. 
\end{bem} 

\section{Application to the topological centre problem} 
We 
first present the proof of Theorem \ref{hau}, for non-compact groups $\G$. -- To establish the 
non-trivial inclusion, let $m \in Z_t(\l1^{**})$. The group $\G$ being non-compact, we infer from 
Lemma \ref{kmazur} that $\l1$ has Mazur's property of level $\kapp(\G)$. So to prove $m \in \l1$, let 
$(h_\al)_{\al \in I} \subseteq \li$ be a bounded net converging $w^*$ to $0$, where $|I|\leq\kapp(\G)$. 
Thanks to Lemma \ref{fakt}, we have the factorization 
$$h_\al = \psi_\al \Diamond h = 
\widetilde{\psi_\al} \odot h \quad (\al \in I)$$ 
with $\psi_\al \in {\overline{\delta_\G}}^{w^*} \subseteq \Ball(\lu)$ and $h \in \li$. -- Here, 
$\widetilde{\psi_\al}$ denotes some arbitrarily chosen Hahn-Banach extension of $\psi_\al$ to $\li^*$. 
We have to show that $a_\al := \langle m, h_\al \rangle \stackrel{\al}{\longrightarrow} 0$. Due to the 
boundedness of $(h_\al)_\al$, it suffices to prove that every convergent subnet of $(a_\al)_\al$ tends to $0$. Let 
$(\langle m, h_{\al_\beta} \rangle)_\beta$ be such a convergent subnet. Furthermore, let 
$$E := \widetilde{w^*-\lim_\gamma \psi_{\al_{\beta_\gamma}}} \in \Ball(\li^*)$$
be an arbitrarily chosen Hahn-Banach extension of some $w^*$-cluster point of the net $(\psi_{\al_{\beta}})_\beta$ 
in $\Ball(\lu)$. 
\par 
We first note that 
$E \odot h=0$, since for arbitrary $g \in \l1$ we obtain: 
\begin{eqnarray*} 
\langle E \odot h, g \rangle &=& \langle E, h \odot g \rangle = 
\lim_\gamma \langle \psi_{\al_{\beta_\gamma}}, h \odot g \rangle \\ 
&=& \lim_\gamma \langle \psi_{\al_{\beta_\gamma}} \Diamond h, g \rangle 
= \lim_\gamma \langle h_{\al_{\beta_\gamma}}, g \rangle \\ 
&=& 0. 
\end{eqnarray*}
Now we conclude that (using twice the assumption $m \in Z_t(\l1^{**})$ and in particular the fact that 
$h . m \in \l1^* \odot \l1 = \LUC(\G)$, cf.\ Lemma 3.1 a) in \cite{laul}) 
\begin{eqnarray*} 
\lim_\beta \langle m, h_{\al_\beta} \rangle &=& 
\lim_\beta \langle m, \widetilde{\psi_{\al_\beta}} \odot h \rangle = 
\lim_\beta \langle m \odot \widetilde{\psi_{\al_\beta}}, h \rangle \\ 
&=& \lim_\beta \langle m . \widetilde{\psi_{\al_\beta}}, h \rangle = 
\lim_\beta \langle \widetilde{\psi_{\al_\beta}}, h . m \rangle \\ 
&=& \lim_\beta \langle \psi_{\al_{\beta}}, h . m \rangle = 
\lim_\gamma \langle \psi_{\al_{\beta_{\gamma}}}, h . m \rangle \\ 
&=& \left\langle w^*-\lim_\gamma \psi_{\al_{\beta_{\gamma}}}, h . m \right\rangle = 
\langle E, h . m \rangle \\ 
&=& \langle m . E, h \rangle = 
\langle m \odot E, h \rangle \\ 
&=& \langle m, E \odot h \rangle = 0 , 
\end{eqnarray*} 
which gives the desired convergence. 
\par 
We now turn to the proof of Theorem \ref{hau10}. 
-- The inclusion $\M(\G) \subseteq Z_t (\lu)$ being 
immediate and classical (cf.\ \cite{wo}, Lemma 3.1), we restrict our attention to the reverse inclusion, where 
in order to avoid trivialities, 
we 
assume the group $\G$ to be non-compact. 
Fix $m \in Z_t (\lu)$. Then, according to Proposition \ref{ftcsatz}, we only have to prove that $\phi (m) \in \B (\li)$ is 
$w^*$-$w^*$-continuous. To this end, 
consider a bounded net $( h_\al )_{\al \in I} \subseteq \li$, where $|I| \leq \kapp (\G)$, 
which converges $w^*$ to $0$. 
Lemma \ref{fakt} 
yields the factorization 
$$h_\al = \psi_\al \Diamond h \quad (\al \in I)$$ 
with $\psi_\al \in {\overline{\delta_\G}}^{w^*} \subseteq \Ball(\lu)$ and $h \in \li$. 
By Lemma \ref{kmazur}, we only have to show that $\phi(m) (h_\al)$ tends $w^*$ to $0$. The latter net 
being bounded, it suffices to prove that every convergent subnet 
tends to $0$. Let 
$(\phi(m) (h_{{\al}_{\beta}}))_\beta$ 
be such a convergent subnet. Furthermore, let 
$$E := w^*-\lim_\gamma \psi_{\al_{\beta_\gamma}} \in {\overline{\delta_\G}}^{w^*} \subseteq \Ball(\lu)$$ 
be a 
$w^*$-cluster point of the net $(\psi_{\al_{\beta}})_\beta$ in $\Ball(\lu)$. 
\par 
We remark that 
$E \Diamond h=0$, since for arbitrary $g \in \l1$ we obtain: 
\begin{eqnarray*} 
\langle E \Diamond h, g \rangle &=& \langle E, h \odot g \rangle = 
\lim_\gamma \langle \psi_{\al_{\beta_\gamma}}, h \odot g \rangle \\ 
&=& \lim_\gamma \langle \psi_{\al_{\beta_\gamma}} \Diamond h, g \rangle 
= \lim_\gamma \langle h_{\al_{\beta_\gamma}}, g \rangle \\ 
&=& 0. 
\end{eqnarray*}
In order to conclude, we will only require the following fact concerning the compatibility 
of our various module operations, which is easy to verify (cf.\ \S 2 in \cite{lapy}): For $\psi \in \lu$, 
$h \in \li$ and $g \in \l1$, we have 
$$(\psi \Diamond h) \odot g = \psi \cdot (h \odot g) .$$ 
Now we obtain: 
\begin{eqnarray*} 
\lim_\beta \langle \phi(m)(h_{\al_{\beta}}) , g \rangle &=& \lim_\beta \langle m \Diamond h_{\al_{\beta}} , g \rangle \\ 
&=& \lim_\beta \langle m, h_{\al_{\beta}} \odot g \rangle \\ 
&=& \lim_\beta \langle m, (\psi_{\al_{\beta}} \Diamond h) \odot g \rangle \\ 
&=& \lim_\beta \langle m, \psi_{\al_{\beta}} \cdot (h \odot g) \rangle \\ 
&=& \lim_\beta \langle m \cdot \psi_{\al_{\beta}}, h \odot g \rangle \\ 
&=& \lim_\gamma \langle m \cdot \psi_{\al_{\beta_\gamma}}, h \odot g \rangle \\ 
&=& \langle m \cdot E, h \odot g \rangle \quad \mathrm{since} ~m \in Z_t(\lu) \\ 
&=& \langle m, E \cdot (h \odot g) \rangle \\ 
&=& \langle m, (E \Diamond h) \odot g \rangle \\ 
&=& 0 , 
\end{eqnarray*} 
which finishes the proof.

\vspace{0.1cm} 
{\sc{Author's address:}} 
\\ 
{\textit{
Department of Mathematical Sciences\\ 
University of Alberta\\ 
Edmonton, Alberta\\ 
Canada T6G 2G1\\ 
E-mail:}} 
mneufang@math.ualberta.ca 
\end{document}